\documentclass{article}

\usepackage{amsmath,amssymb,amsfonts}
\usepackage{amsthm}
\usepackage{graphicx}

\usepackage[margin=2.54cm]{geometry}

\newtheorem{theorem}{Theorem}[section]
\newtheorem{lemma}[theorem]{Lemma}
\newtheorem{proposition}[theorem]{Proposition}

\numberwithin{equation}{section}

\theoremstyle{remark}
\newtheorem{remark}[theorem]{Remark}
\newtheorem{example}[theorem]{Example}

\theoremstyle{definition}
\newtheorem{definition}[theorem]{Definition}

\newcommand{\Ric}{\mathop{\mathrm{Ric}}\nolimits}
\newcommand{\Ad}{\mathop{\mathrm{Ad}}\nolimits}

\newcommand{\tr}{\mathop{\mathrm{tr}}\nolimits}

\newcommand{\rank}{\mathop{\mathrm{rank}}\nolimits}

\def\SO{\mathrm{SO}}
\def\SU{\mathrm{SU}}
\def\SL{\mathrm{SL}}
\def\Sp{\mathrm{Sp}}
\def\Spin{\mathrm{Spin}}
\def\U{\mathrm{U}}
\def\G2{\mathrm{G}_2}
\def\E6{\mathrm{E}_6}
\def\Etwo{\mathrm{E}(2)}
\def\Eoneone{\mathrm{E}(1,1)}
\def\H3{\mathrm{H}_3}

\author{Timothy Buttsworth\thanks{School of Mathematics and Physics, The University of Queensland, St Lucia,~QLD 4072, Australia}~\thanks{Department of Mathematics, Cornell University, Ithaca,~NY 14853, USA} \\
\small{\texttt{t.buttsworth@uq.edu.au}} \and Artem Pulemotov\footnotemark[1]~\thanks{Artem Pulemotov's research is supported under Australian Research Council's Discovery Projects funding scheme (DP180102185).} \\
\small{\texttt{a.pulemotov@uq.edu.au}}}
\title{The prescribed Ricci curvature problem for homogeneous metrics}

\begin{document}

\maketitle

\begin{abstract}
The prescribed Ricci curvature problem consists in finding a Riemannian metric $g$ on a manifold~$M$ such that the Ricci curvature of $g$ equals a given $(0,2)$-tensor field~$T$. We survey the recent progress on this problem in the case where $M$ is a homogeneous space.
\end{abstract}

\section{Introduction}

This paper provides a survey of the recent results on the prescribed Ricci curvature problem with a focus on homogeneous metrics. In Section~\ref{sec_PRC_gen}, we formulate the problem and briefly review its history. In Section~\ref{sec_compact}, we discuss the progress achieved in the framework of compact homogeneous spaces in the papers~\cite{AP16,AP19,MGAPsubm,TBAPYRWZsubm}. Also, we comment on applications to the analysis of the Ricci iteration. Section~\ref{sec_noncomp} raises several open questions and discusses the investigation of the problem on non-compact homogeneous spaces.

\section{The prescribed Ricci curvature problem}\label{sec_PRC_gen}

Consider a smooth manifold $M$ and a symmetric (0,2)-tensor field $T$ on~$M$. The prescribed Ricci curvature problem consists in finding a Riemannian metric $g$ such that
\begin{align}\label{eq_basic_PRC}
\Ric g=T,
\end{align}
where $\Ric$ denotes the Ricci curvature. The investigation of this problem is an important segment of geometric analysis with strong ties to evolution equations. For example, DeTurck's work on~\eqref{eq_basic_PRC} underlay his subsequent discovery of the famous DeTurck trick for the Ricci flow. There is kinship between~\eqref{eq_basic_PRC} and the Einstein equation of general relativity.

Several results regarding the local solvability of~\eqref{eq_basic_PRC} are available in the literature. The first such result appeared in DeTurck's paper~\cite{DDT81}. It states that~\eqref{eq_basic_PRC} has a solution in a neighbourhood of a point $o\in M$ provided that $T$ is nondegenerate at~$o$. Alternative proofs were given in~\cite{AB87,JKweb}. The local solvability of~\eqref{eq_basic_PRC} was further pursued in~\cite{JCDDT94,DDTHG99,AP13a,AP13b}.

Even if $T$ is nondegenerate, one may be unable to find $g$ satisfying~\eqref{eq_basic_PRC} on \emph{all} of $M$. As DeTurck and Koiso demonstrated in~\cite{DDTNK84}, as long as $M$ is compact and $T$ is positive-definite, there is a constant $c_T > 0$ such that $c_T T$ is not the Ricci curvature of any metric; see also~\cite{AB86}. Sufficient conditions for the global solvability of~\eqref{eq_basic_PRC} were obtained by several authors in a range of contexts. Namely, Hamilton in~\cite{RH84}, DeTurck in~\cite{DDT85}, Delano\"e in~\cite{PhD03}, and Delay in a series of papers (see~\cite{ED17b} and references therein) proved the invertability of the Ricci curvature near various ``nice" Einstein metrics. Their results relied on versions of the inverse and implicit function theorems. Pina and his collaborators (see~\cite{RPJDS19} and references therein) studied the solvability of~\eqref{eq_basic_PRC} in a conformal class. The second-named author in~\cite{AP13a,AP13b} and Smith in~\cite{GS16} investigated boundary-value problems for~\eqref{eq_basic_PRC}. DeTurck's and Calvaruso's papers~\cite{DDT83,GC07} provide a sample of the research done in the Lorentzian setting.

Several mathematicians who investigated the prescribed Ricci curvature problem did so in the presence of symmetry. More precisely, they assumed that the metric $g$ and the tensor field $T$ in~\eqref{eq_basic_PRC} were invariant under the action of a group $G$ on the manifold~$M$. Hamilton in~\cite{RH84} considered several different $G$ with $M$ being the 3-dimensional sphere~$\mathbb S^3$. Cao and DeTurck in~\cite{JCDDT94} examined the situation where $G=\SO(d)$ and $M=\mathbb R^d$. The second-named author in~\cite{AP13a,AP13b} assumed $G$ acted on $M$ with cohomogeneity~1.

If the manifold $M$ is closed, instead of trying to prove the existence of a metric $g$ with $\Ric g$ equal to~$T$, one should look for a metric $g$ and a positive number $c$ such that
\begin{align}\label{main_PRC}
\Ric g = cT.
\end{align}
This paradigm was originally proposed by Hamilton in~\cite{RH84} and DeTurck in~\cite{DDT85}. To explain it, we consider the problem of finding a metric on the 2-dimensional sphere $\mathbb S^2$ with prescribed positive-definite Ricci curvature~$T$. According to the Gauss--Bonnet theorem, if such a metric exists, then the volume of $\mathbb S^2$ with respect to $T$ equals $4\pi$. It is relatively easy to prove the converse (see~\cite[Corollary~2.2]{DDT82} and also~\cite[Theorem~2.1]{RH84}). Consequently, it is always possible to find $g$ and $c$ such that~\eqref{main_PRC} holds on~$\mathbb S^2$. In fact, in this setting, $c$ is uniquely determined by~$T$. Hamilton suggests in~\cite[Section~1]{RH84} that the purpose of $c$ is to
compensate for the invariance of the Ricci curvature under scaling of the metric. Note that the shift of focus from~\eqref{eq_basic_PRC} to~\eqref{main_PRC} may be unreasonable on an open manifold or a manifold with non-empty boundary.

For rather detailed surveys of older results on the prescribed Ricci curvature problem, see~\cite[Chapter~5]{AB87} and~\cite[Section~6.5]{TA98}.

\section{Compact homogeneous spaces}\label{sec_compact}

In the paper~\cite{AP16}, the second-named author initiated the systematic investigation of equation~\eqref{main_PRC} for homogeneous metrics. Consider a compact connected Lie group $G$ and a closed connected subgroup $H<G$. Let $M$ be the homogeneous space~$G/H$. It will be convenient for us to assume that the natural action of $G$ on $M$ is effective and that the dimension of $M$ is at least~3. We denote by $\mathcal M$ the set of $G$-invariant Riemannian metrics on~$M$ with the manifold structure as in~\cite[pages~6318--6319]{YNERVS07}. The properties of $\mathcal M$ are discussed in~\cite[Section~4.1]{CB04} in great detail.

Let $T$ be a $G$-invariant tensor field on~$M$. Since $M$ is compact, the prescribed Ricci curvature problem for homogeneous metrics on $M$ consists in finding $g\in\mathcal M$ and $c>0$ that satisfy equation~\eqref{main_PRC}. Our main objective in this section is to state a number of theorems concerning the existence of such $g$ and~$c$. The proofs of most of these theorems rely on the variational interpretation of~\eqref{main_PRC} proposed in~\cite{AP16} and given by Lemma~\ref{lem_var} below.

\subsection{The variational interpretation}\label{subsec_var}

The scalar curvature $S(g)$ of a metric $g\in\mathcal M$ is constant on $M$. Therefore, we may interpret $S(g)$ as the result of applying a functional $S:\mathcal M\to\mathbb R$ to $g\in\mathcal M$. Consider the space 
\begin{align*}
\mathcal M_T=\{g\in\mathcal M\,|\,\tr_gT=1\}
\end{align*}
with the manifold structure inherited from~$\mathcal M$.

\begin{lemma}\label{lem_var}
Given $g\in\mathcal M_T$, formula~(\ref{main_PRC}) holds for some $c\in\mathbb R$ if and only if $g$ is a critical point of the restriction of the functional $S$ to $\mathcal M_T$.
\end{lemma}

In most cases, $S|_{\mathcal M_T}$ is unbounded below; see~\cite[Remark~2.2]{MGAPsubm}. Therefore, one typically searches for maxima, local minima and saddle points of $S|_{\mathcal M_T}$ rather than global minima. For a discussion of the properness of $S|_{\mathcal M_T}$, see~\cite[Remark~2.11]{MGAPsubm}.

Lemma~\ref{lem_var} parallels the well-known variational interpretation of the Einstein equation. Indeed, a metric $g\in\mathcal M$ satisfies
\begin{align}\label{Einstein}
\Ric g = \lambda g
\end{align}
for some $\lambda\in\mathbb R$ if and only it is (up to scaling) a critical point of $S$ on
\begin{align*}
\mathcal M_1=\{g\in\mathcal M\,|\,M~\mbox{has volume~1 with respect to}~g\};
\end{align*}
see, e.g.,~\cite[\S1]{MWWZ86}. The restrictions $S|_{\mathcal M_T}$ and $S|_{\mathcal M_1}$ have substantially different properties. For instance, one can easily deduce from formula~\eqref{sc_curv_formula} below that $S|_{\mathcal M_T}$ is necessarily bounded above if $T$ is positive-definite. This fact is essential to many of the results surveyed here. However, by the same token, $S|_{\mathcal M_1}$ is unbounded above in most situations; see~\cite[Theorem~(2.4)]{MWWZ86} and~\cite[Theorem~1.2]{CB04}.

The contrast in the properties of $S|_{\mathcal M_T}$ and $S|_{\mathcal M_1}$ is hardly surprising. Indeed, despite their apparent similarity, equations~\eqref{main_PRC} and~\eqref{Einstein} are very different in nature. In particular, the former is not diffeomorphism-invariant for general~$T$. Intuitively, one may think of the relationship between~\eqref{main_PRC} and~\eqref{Einstein} as the relationship between the Poisson equation~$\Delta u=f$ and the Helmholtz equation $\Delta u=u$.

It is worth noting that~\eqref{eq_basic_PRC}, as well as~\eqref{main_PRC}, admits a variational interpretation; see~\cite[Theorem~3.1]{RH84}. However, there seem to be no existence results for~\eqref{eq_basic_PRC} to date that arise from this variational interpretation.

\subsection{Maximal isotropy}\label{subsec_max_iso}

If the tensor field $T$ is negative-definite, no metric $g\in\mathcal M$ can satisfy~\eqref{main_PRC} with $c>0$. Indeed, the isometry group of $(M,g)$ would have to be finite by Bochner's theorem (see~\cite[Theorem~1.84]{AB87}), which is clearly impossible. If $T$ is negative-semidefinite and \eqref{main_PRC} holds for some $g\in\mathcal M$ with $c>0$, then $G$ must be abelian, by the same token. In this case, every $G$-invariant metric on $M$ is Ricci-flat. If $T$ has mixed signature, the techniques underlying the results in Sections~\ref{subsec_max_iso}, \ref{subsec_nonmax1} and~\ref{subsec_nonmax2} appear to be ineffective. Thus, we will focus on the situation where $T$ is positive-semidefinite and non-zero. In~\cite{AP16}, the second-named author proved the following result.

\begin{theorem}\label{thm_max_isotr}
Suppose $H$ is a maximal connected Lie subgroup of~$G$. If $T$ is positive-semidefinite and non-zero, then there exists a Riemannian metric $g\in\mathcal M_T$ such that $S(g)\ge S(g')$ for all $g'\in\mathcal M_T$. The Ricci curvature of $g$ coincides with $cT$ for some $c>0$.
\end{theorem}

Clearly, $H$ is a maximal connected Lie subgroup of $G$ if the isotropy representation of $G/H$ is irreducible. In this case, the result is obvious. In~\cite{MWWZ86}, Wang and Ziller use variational methods to show that homogeneous spaces satisfying the conditions of Theorem~\ref{thm_max_isotr} necessarily support Einstein metrics. They discuss numerous examples of such spaces. Theorem~\ref{thm_max_isotr} is similar in spirit to their result. The functionals involved exhibit the same asymptotic behaviour.

Our next goal is to state two sufficient conditions for the solvability of~\eqref{main_PRC} in the case where the maximality assumption of Theorem~\ref{thm_max_isotr} does not hold. In order to do so, we need to introduce structure constants of the homogeneous space $G/H$ and an extension of the scalar curvature functional~$S$.

\subsection{The structure constants}

Denote by $\mathfrak g$ and $\mathfrak h$ the Lie algebras of $G$ and $H$, respectively. Choose a scalar product $Q$ on $\mathfrak g$ induced by a bi-invariant Riemannian metric on~$G$. In what follows, $\oplus$ stands for the $Q$-orthogonal sum. Clearly,
\begin{align*}
\mathfrak g=\mathfrak m\oplus\mathfrak h
\end{align*}
for some $\Ad(H)$-invariant space~$\mathfrak m$. The representation $\Ad(H)|_{\mathfrak m}$ is equivalent to the isotropy representation of $G/H$. We standardly identify $\mathfrak m$ with the tangent space $T_HM$.

Choose a $Q$-orthogonal $\Ad(H)$-invariant decomposition
\begin{align}\label{m_decomp}
\mathfrak m=\mathfrak m_1\oplus\cdots\oplus\mathfrak m_s
\end{align}
such that $\mathfrak m_i\ne\{0\}$ and $\Ad(H)|_{\mathfrak m_i}$ is irreducible for each $i=1,\ldots,s$. The space $\mathfrak m$ may admit more than one decomposition of this form. However, by Schur's lemma, the summands $\mathfrak m_1,\ldots,\mathfrak m_s$ are determined uniquely up to order if $\Ad(H)|_{\mathfrak m_i}$ and $\Ad(H)|_{\mathfrak m_j}$ are inequivalent whenever $i\ne j$. Let $d_i$ be the dimension of~$\mathfrak m_i$. It is easy to show that the number $s$ and the multiset $\{d_1,\ldots,d_s\}$ are independent of the chosen decomposition~\eqref{m_decomp}.

Denote by $B$ the Killing form of $\mathfrak g$. For every $i=1,\ldots,s$, because $\Ad(H)|_{\mathfrak m_i}$ is irreducible, there exists $b_i\ge0$ such that
\begin{align*}
B|_{\mathfrak m_i} = -b_iQ|_{\mathfrak m_i}.
\end{align*}
Given $\Ad(H)$-invariant subspaces $\mathfrak u_1$, $\mathfrak u_2$ and $\mathfrak u_3$ of $\mathfrak m$, define a tensor $\Delta(\mathfrak u_1,\mathfrak u_2,\mathfrak u_3)\in\mathfrak u_1\otimes\mathfrak u_2^*\otimes\mathfrak u_3^*$ by setting
\begin{align*}
\Delta(\mathfrak u_1,\mathfrak u_2,\mathfrak u_3)(X,Y)=\pi_{\mathfrak u_1}[X,Y],\qquad X\in\mathfrak u_2,~Y\in\mathfrak u_3,
\end{align*}
where $\pi_{\mathfrak u_1}$ stands for the $Q$-orthogonal projection onto~$\mathfrak u_1$. Let $\langle \mathfrak u_1\mathfrak u_2\mathfrak u_3\rangle$ be the squared norm of $\Delta(\mathfrak u_1,\mathfrak u_2,\mathfrak u_3)$ with respect to the scalar product on $\mathfrak u_1\otimes\mathfrak u_2^*\otimes\mathfrak u_3^*$ induced by $Q|_{\mathfrak u_1}$, $Q|_{\mathfrak u_2}$ and $Q|_{\mathfrak u_3}$. The fact that $Q$ comes from a bi-invariant metric on $G$ implies
\begin{align*}
\langle\mathfrak u_1\mathfrak u_2\mathfrak u_3\rangle=\langle\mathfrak u_{\rho(1)}\mathfrak u_{\rho(2)}\mathfrak u_{\rho(3)}\rangle
\end{align*}
for any permutation $\rho$ of the set $\{1,2,3\}$. Given $i,j,k\in\{1,\ldots,s\}$, denote
\begin{align*}
\langle\mathfrak m_i\mathfrak m_j\mathfrak m_k\rangle=[ijk].
\end{align*}
The numbers $([ijk])_{i,j,k=1}^s$ are called the \emph{structure constants} of the homogeneous space $M$; cf.~\cite[\S1]{MWWZ86}.

\subsection{The scalar curvature functional and its extension}\label{subsec_sc_curv}

Let us write down a convenient formula for the scalar curvature functional $S$ and define an extension of this functional. In what follows, we implicitly identify $g\in\mathcal M$ with the $\Ad(H)$-invariant scalar product induced by $g$ on $\mathfrak m$ via the identification of $T_HM$ and~$\mathfrak m$. Given a flag $\mathfrak v\subset\mathfrak u\subset\mathfrak m$, denote by $\mathfrak u\ominus\mathfrak v$ the $Q$-orthogonal complement of $\mathfrak v$ in~$\mathfrak u$. The scalar curvature of $g\in\mathcal M$ is given by the equality
\begin{align}\label{sc_curv_basf}
S(g)=-\frac12\tr_gB|_{\mathfrak m}-\frac14
|\Delta(\mathfrak m,\mathfrak m,\mathfrak m)|_g^2,
\end{align}
where $|\cdot|_g$ is the norm induced by~$g$ on $\mathfrak m\otimes\mathfrak m^*\otimes\mathfrak m^*$. If the decomposition~\eqref{m_decomp} is such that
\begin{align}\label{g_form}
g=\sum_{i=1}^sx_i\pi_{\mathfrak m_i}^*Q
\end{align}
for some $x_1,\ldots,x_s>0$, then
\begin{align}\label{sc_curv_formula}
S(g)&=\frac12\sum_{i=1}^s\frac{d_ib_i}{x_i}-\frac14\sum_{i,j,k=1}^s[ijk]\frac{x_k}{x_ix_j}.
\end{align}
For the derivation of formulas~\eqref{sc_curv_basf} and~\eqref{sc_curv_formula}, see, e.g.,~\cite[Chapter~7]{AB87} and~\cite[Lemma~3.2]{APYRsubm}. Given~$g\in\mathcal M$, it is always possible to choose~\eqref{m_decomp} in such a way that~\eqref{g_form} holds. For the proof of this fact, see~\cite[page~180]{MWWZ86}.

Suppose $\mathfrak k$ is a Lie subalgebra of $\mathfrak g$ containing $\mathfrak h$ as a proper subset. It will be convenient for us to denote
\begin{align*}
\mathfrak n=\mathfrak k\ominus\mathfrak h,\qquad \mathfrak l=\mathfrak g\ominus\mathfrak k.
\end{align*}
Let $\mathcal M(\mathfrak k)$ be the space of $\Ad(H)$-invariant scalar products on $\mathfrak n$ equipped with the topology it inherits from the second tensor power of~$\mathfrak n^*$. Our conditions for the solvability of~\eqref{main_PRC} will involve an extension of the functional $S$ to~$\mathcal M(\mathfrak k)$. More precisely, define
\begin{align}\label{hat_S_def}
\hat S(h)=-\frac12\tr_hB|_{\mathfrak n}
-\frac12|\Delta(\mathfrak l,\mathfrak n,\mathfrak l)|_{\rm{mix}}^2
-\frac14
|\Delta(\mathfrak n,\mathfrak n,\mathfrak n)|_h^2,\qquad h\in\mathcal M(\mathfrak k),
\end{align}
where $|\cdot|_{\rm{mix}}$ is the norm induced by $Q|_{\mathfrak l}$ and $h$ on $\mathfrak l\otimes\mathfrak n^*\otimes\mathfrak l^*$ and $|\cdot|_h$ is the norm induced by~$h$ on $\mathfrak n\otimes\mathfrak n^*\otimes\mathfrak n^*$. Intuitively, the first and the third term on the right-hand side are the scalar curvature of $h$ viewed as a degenerate metric on $G/H$. The remaining portion of $\hat S(h)$ quantifies the extent to which $h$ ``notices" the interaction  between $\mathfrak n$ and~$\mathfrak l$. If $\mathfrak k=\mathfrak g$, we identify $\mathcal M(\mathfrak k)$ with $\mathcal M$. In this case, the second term on the right-hand side of~\eqref{hat_S_def} vanishes, and $\hat S(h)$ equals~$S(h)$. If the decomposition~\eqref{m_decomp} is such that
\begin{align}\label{n_h_diag}
\mathfrak k=\Big(\bigoplus_{i\in J_{\mathfrak k}}\mathfrak m_i\Big)\oplus\mathfrak h,\qquad h=\sum_{i\in J_{\mathfrak k}}y_i\pi_{\mathfrak m_i}^*Q,\qquad y_i>0,
\end{align}
for some $J_{\mathfrak k}\subset\{1,\ldots,s\}$, then
\begin{align*}
\hat S(h)&=\frac12\sum_{i\in J_{\mathfrak k}}\frac{d_ib_i}{y_i}-\frac12\sum_{i\in J_{\mathfrak k}}\sum_{j,k\in J_{\mathfrak k}^c}\frac{[ijk]}{y_i}
-\frac14\sum_{i,j,k\in J_{\mathfrak k}}[ijk]\frac{y_k}{y_iy_j},
\end{align*}
where $J_{\mathfrak k}^c$ is the complement of $J_{\mathfrak k}$ in $\{1,\ldots,s\}$; see~\cite[Lemma~2.19]{MGAPsubm}.

\subsection{Non-maximal isotropy: the first existence theorem}\label{subsec_nonmax1}

In order to formulate our first existence theorem for~\eqref{main_PRC} in the case where $H$ is not maximal in $G$, we need a lemma and a definition. As before, let $\mathfrak k$ be a Lie subalgebra of $\mathfrak g$ containing $\mathfrak h$ as a proper subset. Denote 
\begin{align*}
\mathcal M_T(\mathfrak k)=\{h\in\mathcal M(\mathfrak k)\,|\,\tr_hT|_{\mathfrak n}=1\}.
\end{align*} Throughout Sections~\ref{subsec_nonmax1} and~\ref{subsec_nonmax2}, we assume $T$ is positive-definite.

\begin{lemma}\label{prop_sigma_finite}
The quantity $\sigma(\mathfrak k,T)$ defined by the formula
\begin{align*}
\sigma(\mathfrak k,T)=\sup\{\hat S(h)\,|\,h\in\mathcal M_T(\mathfrak k)\}
\end{align*}
satisfies
\begin{align*}
0\le\sigma(\mathfrak k,T)<\infty.
\end{align*}
\end{lemma}

For the proof, see~\cite[Proposition~2.2]{AP19}.

\begin{definition}\label{def_T-apical}
We call $\mathfrak k$ a $T$-\emph{apical} subalgebra of $\mathfrak g$ if $\mathfrak k$ meets the following requirements:
\begin{enumerate}
\item
The inequality $\mathfrak k\ne\mathfrak g$ holds.

\item
There exists a scalar product $h\in\mathcal M_T(\mathfrak k)$ such that 
\begin{align*}
\hat S(h)=\sigma(\mathfrak k,T).
\end{align*}
\item
If $\mathfrak s$ is a maximal Lie subalgebra of $\mathfrak g$ containing $\mathfrak h$, then
\begin{align*}
\sigma(\mathfrak s,T)\le\sigma(\mathfrak k,T).
\end{align*}
\end{enumerate}
\end{definition}
Intuitively, $\mathfrak k$ is $T$-apical if the metrics supported on $\mathfrak k$ have the largest (modified) scalar curvature. The definition also requires that the supremum of $\hat S$ on the set of such metrics be attained. The following result was proven by the second-named author in~\cite{AP19}. It provides a sufficient condition for the solvability of~\eqref{main_PRC}. As before, $\mathfrak l$ stands for~$\mathfrak g\ominus\mathfrak k$.

\begin{theorem}\label{thm_gen}
Suppose that every maximal Lie subalgebra $\mathfrak s$ of $\mathfrak g$ such that $\mathfrak h\subset\mathfrak s$ satisfies the following requirement:
if $\mathfrak u\subset\mathfrak s\ominus\mathfrak h$ and $\mathfrak v\subset\mathfrak g\ominus\mathfrak s$ are non-zero $\Ad(H)$-invariant spaces, then the representations $\Ad(H)|_{\mathfrak u}$ and $\Ad(H)|_{\mathfrak v}$ are inequivalent. Let $\mathfrak k$ be a $T$-apical subalgebra of~$\mathfrak g$ for some positive-definite $T\in\mathcal M$. If
\begin{align}\label{ineq_thm1}
4\sigma(\mathfrak k,T)\tr_Q T|_{\mathfrak l}<-2\tr_QB|_{\mathfrak l}-\langle\mathfrak l\mathfrak l\mathfrak l\rangle,
\end{align}
then there exists $g\in\mathcal M_T$ such that $S(g)\ge S(g')$ for all~$g'\in\mathcal M_T$. The Ricci curvature of~$g$ equals $cT$ for some~$c>0$.
\end{theorem}

As we mentioned in Section~\ref{subsec_var}, the functional $S$ is bounded above on $\mathcal M_T$. It is, therefore, enough to show that its supremum is attained. Estimates imply that, roughly speaking, $S$ drops down to $-\infty$ in all but a finite number of directions. Each such direction is associated with a maximal Lie subalgebra $\mathfrak k$ of~$\mathfrak g$, and the corresponding asymptotic value of $S$ is~$\sigma(\mathfrak k,T)$. If $\mathfrak k$ is $T$-apical, this value is the highest. Formula~\eqref{ineq_thm1} ensures that the slope of $S$ is negative as it approaches~$\sigma(\mathfrak k,T)$. This ensures the existence of a global maximum.

The class of homogeneous spaces satisfying the conditions of Theorem~\ref{thm_gen} is extensive. We discuss examples in Section~\ref{subsec_Examples}.

The above result is moot if $\mathfrak g$ does not have any $T$-apical subalgebras. However, as the second-named author showed in~\cite{AP19}, at least one such subalgebra exists in most situations.

The quantity on the right-hand side of~\eqref{ineq_thm1} is necessarily non-negative; see~\cite[Lemma~2.15]{MGAPsubm}. It is useful to rewrite this quantity in terms of the structure constants of~$M$. Suppose the decomposition~\eqref{m_decomp} is such that the first formula in~\eqref{n_h_diag} holds
for some $J_{\mathfrak k}\subset\{1,\ldots,s\}$. Then
\begin{align*}
-2\tr_QB|_{\mathfrak l}-\langle\mathfrak l\mathfrak l\mathfrak l\rangle=2\sum_{i\in J_{\mathfrak k}^c}d_ib_i-\sum_{i,j,k\in J_{\mathfrak k}^c}[ijk].
\end{align*}

Perhaps, the biggest challenge in applying Theorem~\ref{thm_gen} is to compute the quantity $\sigma(\mathfrak k,T)$ for a given Lie subalgebra~$\mathfrak k$. Our next result  helps overcome this challenge in a number of important cases. For the proof, see~\cite[Proposition~2.5]{AP19}. As above, $\mathfrak n$ stands for~$\mathfrak k\ominus\mathfrak h$.

\begin{proposition}\label{prop_irred}
If $\Ad(H)|_{\mathfrak n}$ is irreducible, then $\mathcal M_T(\mathfrak k)$ consists of a single point. In this case, 
\begin{align}\label{sigma_irred}
\sigma(\mathfrak k,T)=-\frac{2\tr_QB|_{\mathfrak n}+\langle\mathfrak n\mathfrak n\mathfrak n\rangle+2\langle\mathfrak n\mathfrak l\mathfrak l\rangle}{4\tr_QT|_{\mathfrak n}}.
\end{align}
\end{proposition}

\begin{remark}
Suppose $K$ is the connected Lie subgroup of $G$ whose Lie algebra equals~$\mathfrak k$. The irreducibility assumption on the representation $\Ad(H)|_{\mathfrak n}$ in Proposition~\ref{prop_irred} means that the homogeneous space $K/H$ is isotropy irreducible.
\end{remark}

It is useful to restate~\eqref{sigma_irred} in terms of the structure constants of~$M$. If $\Ad(H)|_{\mathfrak n}$ is irreducible, we can choose the decomposition~\eqref{m_decomp} so that $\mathfrak n=\mathfrak m_1$. In this case,
\begin{align*}
T|_{\mathfrak n}=z_1Q|_{\mathfrak m_1},\qquad z_1>0.
\end{align*}
Equality~\eqref{sigma_irred} becomes
\begin{align*}
\sigma(\mathfrak k,T)=\frac1{d_1z_1}\bigg(\frac12d_1b_1-\frac14[111]-\frac12\sum_{j,k\in\{2,\ldots,s\}}[1jk]\bigg).
\end{align*}

\subsection{Non-maximal isotropy: the second existence theorem}\label{subsec_nonmax2}

Our next result, Theorem~\ref{thm_PRC}, provides one more sufficient condition for the solvability of~\eqref{main_PRC}. In some situations, the quantity $\sigma(\mathfrak k,T)$ is difficult to compute for all $\mathfrak k$ even with Proposition~\ref{prop_irred} at hand. Accordingly, it is problematic to determine which $\mathfrak k$ are $T$-apical and to verify condition~\eqref{ineq_thm1}. In such situations, Theorem~\ref{thm_PRC} may be more effective than Theorem~\ref{thm_gen}.

Let us introduce some additional notation and state a definition. Given a bilinear form $R$ on an $\Ad(H)$-invariant non-zero subspace $\mathfrak u\subset\mathfrak m$, define
\begin{align*}
\lambda_-(R)&=\inf\{R(X,X)\,|\,X\in\mathfrak u~\mbox{and}~Q(X,X)=1\},\\
\omega(\mathfrak u)&=\min\{\dim\mathfrak v\,|\,\mathfrak v~\mbox{is a non-zero}\,\Ad(H)\mbox{-invariant subspace of}~\mathfrak u\}.
\end{align*}
Clearly, $\lambda_-(R)$ is the smallest eigenvalue of the matrix of $R$ in a $Q|_{\mathfrak u}$-orthonormal basis of $\mathfrak u$. The number $\omega(\mathfrak u)$ always lies between~1 and~$\dim\mathfrak u$. In fact, $\omega(\mathfrak u)$ is equal to $\dim\mathfrak u$ if $\Ad(H)|_{\mathfrak u}$ is irreducible.

Let $\mathfrak k$ and $\mathfrak k'$ be Lie subalgebras of $\mathfrak g$ such that 
\begin{align}\label{flag}
\mathfrak g\supset\mathfrak k\supset\mathfrak k'\supset\mathfrak h.
\end{align}
Denote
\begin{align*}
\mathfrak l=\mathfrak g\ominus\mathfrak k,\qquad \mathfrak l'=\mathfrak g\ominus\mathfrak k',\qquad \mathfrak j=\mathfrak k\ominus\mathfrak k',\qquad \mathfrak n'=\mathfrak k'\ominus\mathfrak h.
\end{align*}

\begin{definition}\label{def_simple_chain}
We call~(\ref{flag}) a \emph{simple chain} if $\mathfrak k'$ is a maximal Lie subalgebra of $\mathfrak k$ and $\mathfrak h\ne \mathfrak k'$.
\end{definition}

We emphasise that this definition allows the equality $\mathfrak k=\mathfrak g$ but not $\mathfrak k'=\mathfrak k$. In~\cite{MGAPsubm}, Gould and the second-named author proved the following sufficient condition for the solvability of~\eqref{main_PRC}.

\begin{theorem}\label{thm_PRC}
Suppose that every Lie subalgebra $\mathfrak s\subset\mathfrak g$ such that $\mathfrak h\subset\mathfrak s$ and $\mathfrak h\ne\mathfrak s$ meets the following requirements:
\begin{enumerate}
\item
The representations $\Ad(H)|_{\mathfrak u}$ and $\Ad(H)|_{\mathfrak v}$ are inequivalent for every pair of non-zero $\Ad(H)$-invariant spaces $\mathfrak u\subset\mathfrak s\ominus\mathfrak h$ and $\mathfrak v\subset\mathfrak g\ominus\mathfrak s$.

\item
The commutator $[\mathfrak r,\mathfrak s]$ is non-zero for every $\Ad(H)$-invariant 1-dimensional subspace $\mathfrak r$ of $\mathfrak g\ominus\mathfrak s$.
\end{enumerate}
Given a positive-definite $T\in\mathcal M$, if the inequality
\begin{align}\label{ineq_main_thm}
\frac{\lambda_-(T|_{\mathfrak n'})}{\tr_Q T|_{\mathfrak j}}>\frac{2\tr_QB|_{\mathfrak n'}+2\langle\mathfrak n'\mathfrak l'\mathfrak l'\rangle+\langle\mathfrak n'\mathfrak n'\mathfrak n'\rangle}{\omega(\mathfrak n')(2\tr_QB|_{\mathfrak j}+\langle\mathfrak j\mathfrak j\mathfrak j\rangle+2\langle\mathfrak j\mathfrak l\mathfrak l\rangle)}
\end{align}
holds for every simple chain of the form~(\ref{flag}), then there exists a Riemannian metric $g\in\mathcal M_T$ such that $S(g)\ge S(g')$ for all $g'\in\mathcal M_T$. The Ricci curvature of $g$ coincides with $cT$ for some $c>0$.
\end{theorem}

As we mentioned in Section~\ref{subsec_var}, the functional $S$ is bounded above on $\mathcal M_T$. It is, therefore, enough to show that its supremum is attained.

We will discuss examples of homogeneous spaces satisfying the conditions of Theorem~\ref{thm_PRC} in Section~\ref{subsec_Examples}. Under these conditions, the denominator of the fraction on the right-hand side of~\eqref{ineq_main_thm} cannot equal~0, and the fraction itself is necessarily non-negative; see~\cite[Lemma~2.15]{MGAPsubm}. It is possible to rewrite~\eqref{ineq_main_thm} in terms of the structure constants of~$M$, as we rewrote several expressions above; see~\cite[Sections~2.2--2.3]{MGAPsubm}. However, the formulas turn out to be quite bulky, and we will not present them here.

As explained above, both Theorem~\ref{thm_gen} and Theorem~\ref{thm_PRC} provides sufficient conditions for the solvability of~\eqref{main_PRC}. One advantage of the former result over the latter is that it imposes lighter assumptions on the homogeneous space~$M$. In the case where both can be used, Theorem~\ref{thm_gen} seems to yield better conclusions; see~\cite[Remarks~5.2 and~5.5]{AP19}.

\subsection{The case of two isotropy summands}\label{subsec_2sum}

Let us discuss the prescribed Ricci curvature problem on $M$ assuming that the isotropy representation of $M$ splits into two inequivalent irreducible summands. The conditions for the solvability of~\eqref{main_PRC} given by Theorems~\ref{thm_gen} and~\ref{thm_PRC} are not only sufficient but also necessary when this assumption holds. Moreover, the pair $(g,c)\in\mathcal M\times(0,\infty)$ satisfying~\eqref{main_PRC} is unique up to scaling of~$g$. Homogeneous spaces with two irreducible isotropy summands were classified in~\cite{WDMK08,CH12}. Several authors have studied their geometry in detail; see, e.g.,~\cite{AAIC11,MB14,APYRsubm}.

Let $T$ be positive-semidefinite and non-zero. Suppose $s=2$ in every decomposition of the form~\eqref{m_decomp}, i.e.,
\begin{align}\label{m2_dec}
\mathfrak m=\mathfrak m_1\oplus\mathfrak m_2.
\end{align}
According to Theorem~\ref{thm_max_isotr}, equation~\eqref{main_PRC} necessarily has a solution if $\mathfrak h$ is maximal in~$\mathfrak g$. By~\cite[Lemma~4.6]{APYRsubm}, the pair $(g,c)\in\mathcal M\times(0,\infty)$ is unique up to scaling of~$g$. In this section, we assume that there exists a Lie subalgebra $\mathfrak s\subset\mathfrak g$ such that
\begin{align*}
\mathfrak g\supset\mathfrak s\supset\mathfrak h,\qquad \mathfrak h\ne\mathfrak s,\qquad \mathfrak s\ne\mathfrak g.
\end{align*}
Let $\Ad(H)|_{\mathfrak m_1}$ and $\Ad(H)|_{\mathfrak m_2}$ be inequivalent. Without loss of generality, suppose
\begin{align*}
\mathfrak s=\mathfrak m_1\oplus\mathfrak h.
\end{align*}
If $\mathfrak m_2\oplus\mathfrak h$ is also a Lie subalgebra of $\mathfrak g$, then the structure constants $[112]$ and $[221]$ vanish. In this case, all the metrics in $\mathcal M$ have the same Ricci curvature, and the analysis of~\eqref{main_PRC} is easy; see, e.g.,~\cite[Section~4.2]{APYRsubm}. Therefore, we assume $\mathfrak m_2\oplus\mathfrak h$ is not closed under the Lie bracket. This implies $[221]>0$. There exist $z_1,z_2\ge0$ such that $z_1z_2\ne0$ and
\begin{align}\label{T2_def}
T=z_1\pi_{\mathfrak m_1}^*Q+z_2\pi_{\mathfrak m_2}^*Q.
\end{align}

Clearly, if $\mathfrak u\subset\mathfrak s\ominus\mathfrak h$ and $\mathfrak v\subset\mathfrak g\ominus\mathfrak s$ are non-zero $\Ad(H)$-invariant spaces, then $\mathfrak u=\mathfrak m_1$ and $\mathfrak v=\mathfrak m_2$. Since $\Ad(H)|_{\mathfrak m_1}$ and $\Ad(H)|_{\mathfrak m_2}$ are inequivalent, the hypotheses of Theorems~\ref{thm_gen} and~\ref{thm_PRC} are satisfied. Our next result, proven by the second-named author in~\cite{AP16}, settles the prescribed Ricci curvature problem for homogeneous metrics on $M$ under the conditions of this section. Both inequalities~\eqref{ineq_thm1} and~\eqref{ineq_main_thm} reduce to formula~\eqref{2sum_cond} below; cf.~\cite[Section~3]{MGAPsubm}. Thus, Theorems~\ref{thm_gen} and~\ref{thm_PRC} are optimal in the current setting.

\begin{proposition}\label{prop_2sum}
The following statements are equivalent:
\begin{enumerate}
\item
There exist a metric $g\in\mathcal M$ and a number $c>0$ such that the Ricci curvature of $g$ coincides with~$cT$.
\item
The inequality
\begin{align}\label{2sum_cond}
(2b_2d_1d_2-d_1[222])z_1>(2b_1d_1d_2-2d_2[122]-d_2[111])z_2
\end{align}
is satisfied.
\end{enumerate}
When these statements hold, the pair $(g,c)\in\mathcal M\times(0,\infty)$ is unique up to scaling of $g$.
\end{proposition}

The method behind the proof of Proposition~\ref{prop_2sum} in~\cite{AP16} seems to be effective even if $T$ has mixed signature. However, we will not consider this case here.

\subsection{Homogeneous spheres}\label{subsec_spheres}

Exploiting the results in Sections~\ref{subsec_nonmax1} and~\ref{subsec_2sum}, we can obtain a wealth of information about the solvability of~\eqref{main_PRC} for homogeneous metrics on the sphere $\mathbb S^d$. The nature of this information depends strongly on the dimension~$d$. Transitive actions of Lie groups on~$\mathbb S^d$ were classified in the 1940s by Borel, Montgomery and Samelson; see, e.g.,~\cite[Example~6.16]{MARB15}. This classification shows that the homogeneous space $M$ is a sphere if and only if $G$, $H$ and $d$ are as in Table~\ref{table_act}. It will be convenient for us to assume that $Q$ is induced by the round metric of curvature~1.
%We will address the prescribed Ricci curvature problem for metrics invariant under each of the groups $G$ in Table~\ref{table_act}.

\begin{table}
 \caption{Homogeneous structures on spheres. Here, $n\ge2$ in cases~(1)--(2) and $n\ge1$ in cases~(3)--(6).}
\label{table_act}
\centering
  \begin{tabular}{| c | c | c | c |}
    \hline
      & $G$ & $H$ & $d$ \\ \hline
     (1) & $\SO(n+1)$ & $\SO(n)$ & $n$ \\ \hline
     (2) & $\SU(n+1)$ & $\SU(n)$ & $2n+1$ \\ \hline
     (3) & $\U(n+1)$ & $\U(n)$ & $2n+1$ \\ \hline
     (4) & $\Sp(n+1)$ & $\Sp(n)$ & $4n+3$ \\ \hline
     (5) & $\Sp(n+1)\Sp(1)$ & $\Sp(n)\Sp(1)$ & $4n+3$ \\ \hline
     (6) & $\Sp(n+1)\U(1)$ & $\Sp(n)\U(1)$ & $4n+3$ \\ \hline
     (7) & $\SU(2)$ & $\{e\}$ & $3$ \\ \hline     
     (8) & $\SU(2)\SU(2)$ & $\SU(2)$ & $3$ \\ \hline
     (9) & $\Spin(9)$ & $\Spin(7)$ & $15$ \\ \hline
     (10) & $\Spin(7)$ & $\G2$ & $7$ \\ \hline
     (11) & $\G2$ & $\SU(3)$ & $6$ \\ \hline
  \end{tabular}
\end{table}

In cases~(1), (8), (10) and (11), the isotropy representation of $M$ is irreducible. Consequently, all $G$-invariant (0,2)-tensor fields on $M$ are the same up to a scalar multiple, and the analysis of~\eqref{main_PRC} is easy. In cases~(2), (3), (5) and (9), the isotropy representation of $M$ splits into two inequivalent irreducible summands. Given a positive-semidefinite, non-zero, $G$-invariant $T$, formulas~\eqref{m2_dec} and~\eqref{T2_def} hold with $z_1,z_2\ge0$. The spheres $\mathbb S^{2n+1}$, $\mathbb S^{4n+3}$ and $\mathbb S^{15}$ admit the generalised Hopf fibrations
\begin{align*}
\mathbb S^1\hookrightarrow \mathbb S^{2n+1}\to \mathbb {CP}^n, \qquad
\mathbb S^3\hookrightarrow \mathbb S^{4n+3}\to \mathbb{HP}^n,\qquad \mathbb S^7\hookrightarrow \mathbb S^{15}\to
\mathbb S^8.
\end{align*}
We can assume that $\mathfrak m_1$ is vertical and $\mathfrak m_2$ is horizontal. This assumption does not involve any loss of generality. In cases (5) and~(9), Proposition~\ref{prop_2sum} leads to the following result; cf.~\cite[Section~3]{TBAPYRWZsubm}.

\begin{proposition}\label{prop_sph2sum}
Suppose $M$ is as in case~(5) (respectively, case~(9)). Let $T$ be a positive-semidefinite, non-zero and $\Sp(n+1)\Sp(1)$-invariant (respectively, $\Spin(9)$-invariant). A metric $g\in\mathcal M$ solving~(\ref{main_PRC}) for some $c>0$ exists if and only if $T$ satisfies~(\ref{T2_def}) with $(2n+4)z_1>z_2$ (respectively, $14z_1>3z_2$). When it exists, this metric is unique up to scaling.
\end{proposition}

In cases~(2) and~(3), the situation is somewhat different. It is easy to see that every $\SU(n+1)$-invariant metric on $\mathbb S^{2n+1}$ is necessarily $\U(n+1)$-invariant. Proposition~\ref{prop_2sum} leads to the first assertion of Proposition~\ref{prop_SU} below. In case~(7), the space $M$ is a Lie group, and the isotropy representation of $M$ splits into three equivalent 1-dimensional summands. We have the usual Hopf fibration
\begin{align*}
    \mathbb S^1\hookrightarrow \mathbb S^3\to
\mathbb S^2.
\end{align*}
The question of solvability of~\eqref{main_PRC} was settled by Hamilton and the first-named author in~\cite{RH84,TB19}. Specifically, the following result holds.

\begin{proposition}\label{prop_SU}
Let $M$ be a sphere of odd dimension $d\ge3$.
\begin{enumerate}
    \item 
    Assume $M$ is as in case~(2) or~(3). Given a positive-semidefinite, non-zero, $\SU(n+1)$-invariant tensor field~$T$, there exists a metric $g\in\mathcal M$, unique up to scaling, satisfying~(\ref{main_PRC}) for some $c>0$.
    \item
    Assume $M$ is as in case~(7) and $T$ is a positive-semidefinite, non-zero, left-invariant tensor field. If $T$ is non-degenerate, then there is a metric $g\in\mathcal M$, unique up to scaling, satisfying~(\ref{main_PRC}) for some $c>0$. If $\rank T=1$, then there is a 2-parameter family of such metrics. If $\rank T=2$, none exists.
\end{enumerate}
\end{proposition}

\begin{remark}\label{rem_Hopf}
It is possible to state and prove Proposition~\ref{prop_sph2sum} and the first assertion of Proposition~\ref{prop_SU} avoiding the homogeneous space formalism. Instead, one may use the language of Hopf fibrations as in~\cite[Remark~3.2]{TBAPYRWZsubm}.
\end{remark}

Finally, we turn to cases~(4) and~(6). If $M$ equals $\Sp(n+1)/\Sp(n)$, then the isotropy representation of $M$ splits into 4 irreducible summands, i.e.,
\begin{align}\label{m4_dec}
\mathfrak m=\mathfrak m_1\oplus\mathfrak m_2\oplus\mathfrak m_3\oplus\mathfrak m_4.
\end{align}
Three of these summands, say, $\mathfrak m_1$, $\mathfrak m_2$ and $\mathfrak m_3$, are 1-dimensional and equivalent. Given a positive-definite, $\Sp(n+1)$-invariant tensor field $T$, we can choose the decomposition~\eqref{m4_dec} so that
\begin{align}\label{T4_form}
T=z_1\pi_{\mathfrak m_1}^*Q+z_2\pi_{\mathfrak m_2}^*Q+z_3\pi_{\mathfrak m_3}^*Q+q\pi_{\mathfrak m_4}^*Q
\end{align}
for some $z_1,z_2,z_3,q>0$. The following sufficient condition for the solvability of~\eqref{main_PRC} was obtained by the authors, Rubinstein and Ziller in~\cite{TBAPYRWZsubm}.

\begin{theorem}\label{thm_S4n3}
Let $M$ be as in case~(4). Suppose $T$ is a positive-definite, $\Sp(n+1)$-invariant tensor field on $M$ satisfying~(\ref{T4_form}). If
\begin{align}\label{cond_S4n3}
(2n+4)\min\{z_1,z_2,z_3\}>q,    
\end{align}
then there exists a metric $g\in\mathcal M$ such that~(\ref{main_PRC}) holds for some $c>0$.
\end{theorem}

This result was proven in~\cite{TBAPYRWZsubm} by means of topological degree theory. It can also be derived from Theorem~\ref{thm_gen} above. The homogeneous sphere $\Sp(n+1)/\Sp(n)$ is, thus, an example of a space with equivalent isotropy summands on which this theorem applies.

It is possible to replace~\eqref{cond_S4n3} will a less restrictive formula that is more difficult to verify. We will not discuss this here; see~\cite[Theorem~5.4]{TBAPYRWZsubm}.

If $M$ is as in case~(6), then a positive-definite $\Sp(n+1)\U(1)$-tensor field $T$ on $M$ can be written in the form~\eqref{T4_form} with two of the $z_i$s equal to each other. In this situation, again, Theorem~\ref{thm_S4n3} provides a sufficient condition for the solvability of~\eqref{main_PRC}.

The results stated above cover all homogeneous structures on spheres. By analogy with Proposition~\ref{prop_sph2sum} and Remark~\ref{rem_Hopf}, we can use Proposition~\ref{prop_2sum} and generalised Hopf fibrations to address the prescribed Ricci curvature problem on projective spaces. We will not discuss the details here; see~\cite{TBAPYRWZsubm}.

\subsection{Further examples}\label{subsec_Examples}

Aside from homogeneous spheres and projective spaces, Theorems~\ref{thm_gen} and~\ref{thm_PRC} provide easy-to-verify sufficient conditions for the solvability of~\eqref{main_PRC} on large classes of generalised Wallach spaces and generalised flag manifolds. Previous literature contains little information concerning the prescribed Ricci curvature problem on such spaces. However, several other aspects of their geometry have been investigated thoroughly; see the survey~\cite{AA15}.

Theorem~\ref{thm_gen} is an effective tool for solving~\eqref{main_PRC} on generalised Wallach spaces with inequivalent isotropy summands. The reader will find classifications of such spaces in~\cite{ZCYKKL16,YN16}. There are several infinite families and 10 isolated instances, excluding products. We will illustrate the usage of Theorem~\ref{thm_gen} by considering an example. For a more detailed discussion, see~\cite{AP19}.

\begin{example}
Let $M$ be the generalised Wallach space $\E6/(\Sp(3)\Sp(1))$ with the decomposition~\eqref{m_decomp} chosen as in~\cite{YN16}. Suppose that $Q=-B$ and that $T$ is an element of $\mathcal M$. The formula
\begin{align}\label{T_form3}
T=-z_1\pi_{\mathfrak m_1}^*B-z_2\pi_{\mathfrak m_2}^*B-z_3\pi_{\mathfrak m_3}^*B
\end{align}
holds for some $z_1,z_2,z_3>0$.
Theorem~\ref{thm_gen} implies that a Riemannian metric $g\in\mathcal M$ satisfying~\eqref{main_PRC} for some $c>0$ exists if the triple $(z_1,z_2,z_3)$ lies in the set
\begin{align*}
\{(x,y,z)&\in(0,\infty)^3\,|\,3x\le2y,~5x\le6z,~7y+3z<20x\} \\
&\phantom{\in}~\cup\{(x,y,z)\in(0,\infty)^3\,|\,3x\ge2y,~5y\le9z,~21x+18z<52y\} \\
&\phantom{\in}~\cup\{(x,y,z)\in(0,\infty)^3\,|\,5x\ge6z,~5y\ge9z,~5x+10y<36z\}.
\end{align*}
\end{example}

Generalised flag manifolds form an important class of homogeneous spaces with applications across a range of fields. The reader will find classifications in~\cite{AAIC10,SAIC11,AAICYS13,ICYS14}. There are numerous infinite families and isolated instances. As shown by Gould and the second-named author in~\cite{MGAPsubm}, if $M$ is a generalised flag manifold, then $M$ satisfies the hypotheses of Theorem~\ref{thm_PRC}. The simple chains are easy to identify, and inequality~\eqref{ineq_main_thm} is straightforward to check, as long as the structure constants of $M$ are known. For the computation of these constants for a variety of spaces, see~\cite{MK90,AAIC10,AAIC11,SAIC11,AAICYS13,ICYS14}.

Evidently, generalised flag manifolds satisfy the hypotheses of Theorem~\ref{thm_gen} as well. The $T$-apical subalgebras are relatively easy to find, and the verification of~\eqref{ineq_thm1} is straightforward, if the isotropy representation of $M$ splits into~5 or fewer irreducible summands. Let us provide an example. We refer to~\cite[Section~5]{AP19} for more. 

\begin{example}\label{example_G2U2}
Suppose $M$ is the generalised flag manifold $\G2/\U(2)$ in which $\U(2)$ corresponds to the long root of~$\G2$. Let the decomposition~\eqref{m_decomp} be as in~\cite{SAIC11}. Given $T\in\mathcal M$, formula~\eqref{T_form3} holds for some $z_1,z_2,z_3>0$. Theorem~\ref{thm_gen} implies that a Riemannian metric $g\in\mathcal M$ satisfying~\eqref{main_PRC} for some $c>0$ exists if $(z_1,z_2,z_3)$ is in the set
\begin{align*}
\{(x,y,z)&\in(0,\infty)^3\,|\,9y\le2z,~x+z<12y\}\cup\{(x,y,z)\in(0,\infty)^3\,|\, 9y\ge2z,~6x+3y<10z\}.
\end{align*}
\end{example}

\subsection{Ricci iterations}

The results discussed above lead to new existence theorems for Ricci iterations. More precisely, let $(g_i)_{i=1}^\infty$ be  a sequence of Riemannian metrics on a smooth manifold. One calls $(g_i)_{i=1}^\infty$ a \emph{Ricci iteration} if
\begin{align}\label{eq_iter}
\Ric g_{i+1}=g_i
\end{align}
for $i\in\mathbb N$. First considered by Rubinstein in~\cite{YR07}, sequences satisfying~\eqref{eq_iter} have been investigated intensively in the framework of K\"ahler geometry; see the survey~\cite{YR14}, and also the paper~\cite{TDYR19} for the most recent and comprehensive results. One may interpret~\eqref{eq_iter} as a discretisation of the Ricci flow, as explained in, e.g.,~\cite[Section~6]{YR14}. This observation establishes an intriguing link between prescribed curvature problems and the theory of geometric evolutions.

In~\cite{APYRsubm}, the second-named author and Rubinstein initiated the study of Ricci iterations in the homogeneous framework. Relying on Theorem~\ref{thm_max_isotr} and Proposition~\ref{prop_2sum}, they achieved a deep understanding of~\eqref{eq_iter} on spaces with two inequivalent irreducible isotropy summands and proved existence and compactness results for~\eqref{eq_iter} on spaces with maximal isotropy. They also observed an interesting connection between Ricci iterations and ancient solutions to the Ricci flow. On the basis of this work, the authors, Rubinstein and Ziller investigated~\eqref{eq_iter} on homogeneous spheres and projective spaces in~\cite{TBAPYRWZsubm}. One of the results of~\cite{TBAPYRWZsubm} establishes the existence and the convergence of Ricci iterations by means of the inverse function theorem. This approach turns out to be effective even on spaces without a homogeneous structure, as the first-named author and Hallgren demonstrate in~\cite{TBMHforth}.

\section{Open questions and non-compact homogeneous spaces}\label{sec_noncomp}

While the results surveyed in Section~\ref{sec_compact} provide a large amount of information about the prescribed curvature problem for homogeneous metrics, and even settle the problem in some cases, a number of questions remain open. For instance, what necessary conditions for the solvability of~\eqref{main_PRC} can one state on compact spaces with more than two isotropy summands? In what circumstances is it possible to establish the uniqueness of $g$ on such spaces (up to scaling)? What approaches can one use to study the solvability of~\eqref{main_PRC} when $M$ is compact and $T$ has mixed signature? These questions are, in fact, still open on spheres of dimension $4n+3$ acted upon by $\Sp(n+1)$.

\subsection{The non-compact case}

The investigation of the prescribed Ricci curvature problem for homogeneous metrics on \emph{non-compact} spaces is a promising area of research. So far, progress has been scarce. On the other hand, while equations~\eqref{eq_basic_PRC} and~\eqref{main_PRC} appear to be difficult to solve in the homogeneous non-compact setting, a number of related questions are more tractable. For instance, what are the possible signatures of the Ricci curvature of a left-invariant metric on a connected Lie group~$\mathcal G$? In~\cite{Milnor76}, Milnor proposed an answer to this question assuming $\mathcal G$ was 3-dimensional and unimodular. His arguments relied on the observation that any left-invariant metric can be diagonalised in a basis satisfying ``nice" commutation relations. Subsequently, in~\cite{HaLee09}, Ha and Lee obtained counterparts of Milnor's results for non-unimodular~$\mathcal G$. Possible signatures of the Ricci curvature were further studied by Kremlev and Nikonorov on 4-dimensional Lie groups and by Boucetta, Djiadeu Ngaha and Wouafo Kamga on nilpotent groups of all dimensions; see~\cite{KremNik08,BKN16}. Several other researchers addressed questions related to the solvability of~\eqref{eq_basic_PRC} and~\eqref{main_PRC} in indirect ways. Namely, in~\cite{KowNik96}, Kowalski and Nikcevic extended the results of~\cite{Milnor76,HaLee09} by finding all possible triples of numbers that could be eigenvalues of the Ricci curvature of a left-invariant metric in 3 dimensions. In \cite{E08}, Eberlein investigated the space of metrics with prescribed Ricci operator on 2-step nilpotent Lie groups.

Let us make a few observations concerning the solvability of~\eqref{eq_basic_PRC} on non-compact homogeneous spaces. These observations mainly deal with the non-existence of metrics with non-negative Ricci curvature.
%In Section~\ref{subsec_3d_unim}, we discuss the solvability of~\eqref{eq_basic_PRC} on unimodular Lie groups of dimension~3. The result we present extends those of Milnor in~\cite{Milnor76}.
%\subsection{Non-existence of non-negatively curved metrics}\label{subsec_nonex_noncom}
Equality~\eqref{eq_basic_PRC} cannot hold for the metric on a non-compact Riemannian homogeneous space if the tensor field $T$ is positive-definite. Indeed, the Bonnet--Myers theorem implies the following result.

\begin{proposition}\label{MyerTheorem}
If a connected Riemannian homogeneous space has positive-definite Ricci curvature, then it must be compact.
\end{proposition}

Using the Cheeger-Gromoll splitting theorem, one can prove the following extension of this result (see~\cite[Theorem 6.65 and Remark~6.66\,(f)]{AB87}).

\begin{proposition}\label{NNRC}
If a non-compact, connected Riemannian homogeneous space has positive-semidefinite Ricci curvature, then it is isometric to the Riemannian product of a compact Riemannian homogeneous space with positive-semidefinite Ricci curvature and a Euclidean space with standard metric.
\end{proposition}

A natural question arises: can the metric on a non-compact Riemannian homogeneous space satisfy~\eqref{eq_basic_PRC} with $T=0$? The following result, known as the Alekseevskii--Kimel'fel'd theorem, provides the answer.

\begin{theorem}\label{RFIF}
A Ricci-flat Riemannian homogeneous space is isometric to the Riemannian product of a flat torus and a Euclidean space with standard metric.
\end{theorem}

One can derive this result from Proposition~\ref{NNRC} and the Bochner theorem (see~\cite[Theorem~1.84]{AB87}). It was originally proven in~\cite{AK75} by other methods.

Theorem~\ref{RFIF} would be false, of course, without the homogeneity assumption. One can see this, for example, by considering Calabi--Yau manifolds.

\subsection{Unimodular Lie groups of dimension~3}\label{subsec_3d_unim}

In~\cite{TB19}, the first-named author settled the question of the solvability of equation~\eqref{main_PRC} for left-invariant Riemannian metrics on unimodular Lie groups in 3 dimensions. These results provide insight into the prescribed Ricci curvature problem on non-compact homogeneous spaces. We summarise them in Theorem~\ref{thm_unim3} below.

Consider a 3-dimensional connected unimodular Lie group $\mathcal G$ with Lie algebra~$\mathfrak x$. It is well-known that $\mathcal G$ necessarily possesses a Milnor frame, i.e., is a basis $\{V_1,V_2,V_3\}$ of $\mathfrak x$ such that 
\begin{align*}
    [V_i,V_j]=\epsilon_{ijk}\lambda_kV_k,\qquad i,j\in \{1,2,3\},
\end{align*}
where $\epsilon_{ijk}$ is the Levi-Civita symbol, $k\in\{1,2,3\}$ is distinct from $i$ and~$j$, and $\lambda_k\in\mathbb R$. By scaling and reordering, we can assume the following:
\begin{itemize}
\item[(a)]
Each $\lambda_k$ lies in $\{-2,0,2\}$. 
\item[(b)]
There are more non-negative $\lambda_k$s than negative ones.
\item[(c)]
The equality $\lambda_3=0$ holds if any of the $\lambda_k$s is~0.
\item[(d)]
The equality $\lambda_1 = 2$ holds unless $\lambda_k=0$ for all $k$. 
\end{itemize}
With these additional constraints, the triple $(\lambda_1,\lambda_2,\lambda_3)$ is the same for all Milnor frames on $\mathcal G$. Moreover, this triple determines the Lie algebra $\mathfrak x$ uniquely; see, e.g.,~\cite{Milnor76} and the first column in Table~\ref{table_unim}.

Let $T$ be a left-invariant tensor field on $\mathcal G$. In~\cite{TB19}, the first-named author proved the following result.

\begin{theorem}\label{thm_unim3}
A left-invariant Riemannian metric $g$ satisfying~(\ref{main_PRC}) for some $c>0$ exists if and only if one can find a Milnor frame $\{V_1,V_2,V_3\}$ such that the above constraints (a)--(d) and the following three statements hold:
\begin{enumerate}
    \item The tensor field $T$ is diagonalisable in~$\{V_1,V_2,V_3\}$.
    \item The signs of the eigenvalues $(T_1,T_2,T_3)$ of the matrix of $T$ in the basis $\{V_1,V_2,V_3\}$ appear in the second column of Table~\ref{table_unim}.
    \item The eigenvalues $(T_1,T_2,T_3)$ satisfy the conditions in the third column of Table~\ref{table_unim}.
\end{enumerate}
\end{theorem}

It is natural to ask whether the constant $c$ in~\eqref{main_PRC} is uniquely determined by $T$ and whether the metric $g$ satisfying~\eqref{main_PRC} for given $T$ and $c>0$ is in any sense unique. The answers to these questions are provided in the fourth and fifth columns of Table~\ref{table_unim}. The symbol $\sim$ denotes equality up to scaling.

Results of~\cite{Milnor76,HaLee09} on the possible signatures of the Ricci curvature yield the information contained in the first two columns. Thus, Theorem~\ref{thm_unim3} extends these results.

\begin{table}
\caption{Existence and uniqueness of solutions to~\eqref{main_PRC} on $\mathcal G$.}
\label{table_unim}
\centering
\begin{tabular}{ |c|c|c|c|c| }
 \hline
Lie Group  & Signs of & Conditions on 
& Is $c$ the & 
$g_1 \sim g_2$ for any \tabularnewline
$(\lambda_{1},\lambda_{2},\lambda_{3})$ & $(T_1,T_2,T_3)$ & $(T_{1},T_{2},T_{3})$
&same for all & solutions $(g_1,c)$\tabularnewline
& & & solutions? &  and $(g_2,c)$ \tabularnewline
\hline 
$\SO(3)$ or $\SU(2)$ & $(+,+,+)$ & - & Yes & Yes\tabularnewline
$(2,2,2)$ & $(+,0,0)$ & - & Yes & No\tabularnewline
 & $(0,+,0)$ & - & Yes & No\tabularnewline
 & $(0,0,+)$ & - & Yes & No\tabularnewline
  & $(+,-,-)$ & $(T_1+T_2+T_3)^3\ge27T_1T_2T_3$ & No & Yes\tabularnewline
   & $(-,+,-)$ & $(T_1+T_2+T_3)^3\ge27T_1T_2T_3$ & No & Yes\tabularnewline
    & $(-,-,+)$ & $(T_1+T_2+T_3)^3\ge27T_1T_2T_3$ & No & Yes\tabularnewline
\hline 
$\SL(2,\mathbb R)$ & $(+,-,-)$ & $T_{3}+T_{1}>0$ & Yes & Yes\tabularnewline
$(2,2,-2)$ & $(-,+,-)$ & $T_{3}+T_{2}>0$ & Yes & Yes\tabularnewline
 & $(-,-,+)$ & $\max\{ -T_{1},-T_{2}\} <T_{3}~\mbox{or}$ & Yes & Yes\tabularnewline
 &  & $\min\{ -T_{1},-T_{2}\} >T_{3}~\mbox{or}$ & Yes & Yes\tabularnewline
 &  & $T_{3}=-T_{1}=-T_{2}$ & Yes & No\tabularnewline
 & $(-,0,0)$ & - & Yes & No\tabularnewline
 & $(0,-,0)$ & - & Yes & No\tabularnewline
\hline 
$\Etwo$ & $(0,0,0)$ & - & No & No\tabularnewline
$(2,2,0)$ & $(+,-,-)$ & $T_{1}+T_{2}>0$ & Yes & Yes\tabularnewline
 & $(-,+,-)$ & $T_{1}+T_{2}>0$ & Yes & Yes\tabularnewline
\hline 
$\Eoneone$ & $(0,0,-)$ & - & Yes & No\tabularnewline
$(2,-2,0)$ & $(+,-,-)$ & $T_{1}+T_{2}>0$ & Yes & Yes\tabularnewline
 & $(-,+,-)$ & $T_{1}+T_{2}>0$ & Yes & Yes\tabularnewline
\hline 
$\H3$ & $(+,-,-)$ & - & Yes & Yes\tabularnewline
$(2,0,0)$ & & & & \tabularnewline
\hline 
$\mathbb{R}^{3}$ & $(0,0,0)$ & - & No & No\tabularnewline
$(0,0,0)$ & & & & \tabularnewline
\hline
\end{tabular}
\end{table}

\section*{Acknowledgements}

The authors express their gratitude to Jorge Lauret for pointing out a small error in the first version of the manuscript.

\end{document}